\documentclass[12pt]{article}
\input{amssym.def}
\input{amssym}
\usepackage{eufrak}

\newtheorem{theorem}{Theorem}

\newtheorem{proposition}[theorem]{Proposition}
\newtheorem{remark}[theorem]{Remark}

\def\qq{q^{-1}}
\def\pa{\partial}

\def\De{\Delta}

\def\tpa{\tilde{\pa}_t}

\def\hh{{\nu}}

\def\de{\delta}

\def\tD{\tilde{D}}
\def\ot{\otimes}

\def\C{{\Bbb C}}

\def\Sym{{\rm Sym\, }}

\def\vv{V^{\otimes 2}}

\def\lhq{\ifmmode {\cal L}(q,\hbar)\else ${\cal L}(q,\hbar)$\fi}

\def\lqh{\ifmmode {\cal L}(q,\hbar)\else ${\cal L}(q,\hbar)$\fi}

\def\xx{\mathbf{x}}

\def\be{\begin{equation}}
\def\ee{\end{equation}}

\begin{document}

\makeatletter
\renewcommand{\theequation}{{\thesection}.{\arabic{equation}}}
\@addtoreset{equation}{section} \makeatother

\title{Doubles of associative algebras and their applications}

\author{\rule{0pt}{7mm} Dimitri
Gurevich\thanks{gurevich@ihes.fr}\\
{\small\it LMI, UPHF,
59313 Valenciennes, France}\\
{\small \it and}\\
{\small \it
Interdisciplinary Scientific Center J.-V.Poncelet}\\
{\small \it Moscow 119002, Russian Federation}\\
\rule{0pt}{7mm} Pavel Saponov\thanks{Pavel.Saponov@ihep.ru}\\
{\small\it
National Research University Higher School of Economics,}\\
{\small\it 20 Myasnitskaya Ulitsa, Moscow 101000, Russian Federation}\\
{\small \it and}\\
{\small \it
Institute for High Energy Physics, NRC "Kurchatov Institute"}\\
{\small \it Protvino 142281, Russian Federation}}

\maketitle

\begin{abstract}
For a couple of associative algebras we define the notion of their double and give a set of examples. Also, we discuss applications of such doubles to representation theory of
certain quantum algebras and to a new type of Noncommutative  Geometry.
\end{abstract}

\section{Introduction}

In this letter by a double of associative algebras we mean an ordered couple $(A,B)$ of associative unital algebras $A$ and $B$ such that their tensor products $B\ot A$ can be also endowed with
an associative product by means of a {\em permutation map} $\sigma: A\ot B\to B\ot A$. If the algebra $A$ is equipped with a counit (an algebra homomorphism) $\varepsilon:A\to \C$, then under
some natural conditions on $\sigma$ and $\varepsilon$ the algebra $A$ can be represented in the algebra $B$.

The simplest example of such a double is a Heisenberg-Weyl (HW) algebra. The smash-product of a bi-algebra $A$ and an $A$-module $M$ is another example of a double. In this case  the role
of the algebra $B$ can be played by the free tensor algebra $T(M)=\oplus_k M^{\ot k}$ or by some of its quotient algebras. We are mainly interested in doubles related to  braidings.

Let $V$ be a finite dimensional complex vector space, $\dim V = N$. An invertible operator $R:\vv\to\vv$ is called a {\em braiding}, if it is subject to the {\em braid relation}
$$
R_{12}\, R_{23}\,R_{12}\,=R_{23}\,R_{12}\,R_{23},\qquad R_{12}=R\ot I, \quad R_{23}=I\ot R.
$$
Hereafter, $I$ stands for the identity operator or its matrix. A braiding $R$ is called respectively an {\it involutive} or a {\it Hecke symmetry} if it is subject to a supplementary condition
$$
R^2=I\qquad {\rm or}\qquad (qI-R)(\qq I+R)=0 \quad q\not\in\{0,\pm 1\}.
$$
The best known examples of Hecke symmetries come from the Drinfeld-Jimbo Quantum Groups (QG) $U_q(sl(N))$. They are deformations of the usual flips.
Nevertheless, there exist involutive and Hecke symmetries, which are neither deformations of the flips nor super-flips.

All symmetries, we are dealing with, are assumed to be skew-invertible (see \cite{1}). We mainly deal with Hecke symmetries $R=R(q)$ at a {\it generic} value of the parameter $q$.
To any such a Hecke symmetry $R$  we associate $R$-analogs of the symmetric  and skew-sym\-met\-ric  algebras of the space
$V$ by respectively setting
$$
{\mathrm{Sym}}_R(V)=T(V)/\langle {\rm Im}(qI-R)\rangle,\qquad {\mathsf{\Lambda}}_R(V)=T(V)/\langle {\rm Im}(\qq I+R)\rangle.
$$
Besides, we consider the so-called RTT and Reflection Equation (RE) algebras defined respectively by
\be
R_{12} T_{1} T_{2}-T_{1}T_{2} R_{12}=0,
\label{RTT}
\ee
\be R_{12}L_{1}\,R_{12}L_{1}-L_{1}R_{12}L_{1} R_{12}=0,
\label{RE}
\ee
where $T=\|t_i^j\|_{1\leq i,j\leq N}$, $L=\|l_i^j\|_{1\leq i,j\leq N}$, $T_1=T\ot I$ and $T_2=I\ot T$ etc.

In section 3 we exhibit examples of doubles $(A,B)$, where the algebra $A$ is an RTT or an RE algebra. In a number of papers there were considered doubles with  RE algebras
playing the role of $A$ and the corresponding RTT algebras, playing the role of $B$. By contrast, in \cite{2, 3} we considered doubles, where $B$ was another copy of the RE algebra.
Combining the generating matrices of the algebras $A$ and $B$  we constructed other doubles, giving rise to the notion of partial derivatives in the noncommutative generators of $B$.

Since the RE algebra in its modified form tends to the algebra $U(gl(N))$, provided the Hecke symmetry  $R=R(q)$ tends to the usual flip $P$ as $q\to 1$, we obtained  partial derivatives
on the algebra $U(gl(N))$\footnote{Sometimes we will deal with the algebra $U(gl(N)_h)$, where $h$ is a numerical multiplier introduced in the bracket of the Lie algebra $gl(N)$. This rescaling
of the bracket enables us to treat the algebra  $U(gl(N)_h)$ as a quantization of the commutative algebra $\Sym(gl(N))$ with respect to the linear Poisson bracket.}. These partial derivatives
gives rise to a new noncommutative (NC) differential calculus which is $GL(N)$-covariant and which turns into the usual calculus on the algebra $\Sym(gl(N))$ as $h\to 0$.

In a particular case $N=2$ we treat  the compact form $U(u(2)_h)$ of the algebra $U(gl(2)_h)$ as an NC version of the polynomial algebra on the Minkowski space. Given a differential operator
with polynomial coefficients on the classical Minkowski space, we quantize the coefficients and replace the usual partial derivatives with their "quantum counterparts". In this way we get an operator,
defined on the algebra $U(u(2)_h)$ which turns into the initial one as $h\to 0$. We call this procedure {\em the quantization with an NC configuration space}. In \cite{4, 5} we extend this procedure on
some operators with non-polynomial coefficients.

In  the present letter we reproduce some elements of this NC calculus. However, the main our objective is comparing the doubles $(A, B)$ where $A$ is an RTT or RE algebra, and constructing
the corresponding representations of these algebras. In the last section we consider an example of the mentioned quantization with an NC configuration space.

 {\bf Acknowledgement}
The work of P.S.  was partially funded the RFBR grant 19-01-00726.

\section{Representations via doubles, first examples}

Let $A$ and $B$ be two associative unital algebras endowed with a linear map
$ \sigma: A\ot B\to B\ot  A$, such that
$$
\sigma\circ(\mu_{A}\otimes\mathrm{id}_B) =\mu_{A}\circ\sigma_{12}\circ\sigma_{23} \quad{\rm on}\quad A\ot A\ot B,
$$
$$
\sigma\circ(\mathrm{id}_A\otimes\mu_{B}) =\mu_{B}\circ\sigma_{23}\circ\sigma_{12} \quad{\rm on}\quad A\ot B\ot B,
$$
$$
\sigma(1_A\ot b)=b\ot 1_A,\quad \sigma(a\ot 1_B)=1_B\ot a \qquad\forall\, a\in A,\, \forall\,b\in B,
$$
where $\mu_A: A\ot A\to A $ is the product in the algebra $A$, $1_A$ is its unit, and similarly for $B$.

Under these assumptions the space $B\ot A$ can be equipped with a bilinear map $*$:
$$
(B\ot A)^{\ot 2}\stackrel{*}{\rightarrow} B\ot A:\quad
(b\ot a) * (b'\ot a'):=(\mu_B\otimes \mu_A)\circ (\mathrm{id}_B\otimes\sigma_{23}\otimes\mathrm{id}_A) (b\ot a\ot b'\ot a').
$$

\begin{proposition}
The  map $*$ endows the space $B\otimes A$ with the structure of a unital associative algebra with the unit element $1_B\ot 1_A$.
\end{proposition}
We call the corresponding algebra {\it the double of associative algebras} $A$ and $B$ and denote it as $B\ot_{\sigma} A$.

If the algebra $A$ is equipped with a counit (an algebra homomorphism) $\varepsilon_A:A\to \C$, then we define an action of the algebra $A$ onto  $B$ by
the rule
$$
a\triangleright b = (\mathrm{id}_B\otimes\varepsilon_A)\circ \sigma(a\ot b),\quad \forall\,a\in A,\, \forall\,b \in B.
$$
Identifying $b\ot 1_\C$ and $b$, we get that $a\triangleright b\in B$, so each element $a\in A$ defines a linear operator
$$
Op(a):B\to B.
$$

\begin{proposition}
The map $a\mapsto Op(a)$ defines a representation of the algebra $A$ in the algebra $B$:
$$
Op(ab)=Op(a) Op(b),\quad Op(1_A)=I.
$$
\end{proposition}
Note that if $\sigma=P$, this representation becomes trivial  $a  \mapsto \varepsilon(a)\,I$.

\medskip

\noindent
1. As an example we consider an  HW algebra, generated by two polynomial subalgebras $A=\C[x^1, ..., x^m]$ and $B=\C[x_1, ..., x_m]$.
 Introduce the following permutation relations\footnote{By permutation relations we mean equalities $a\ot b=\sigma(a\ot b)$, $a\in A$, $b \in B$. All the doubles $(A,B)$ below
are defined via relations on  generators of each component  and the permutation relations.}
\be
x^j\, x_i=x_i\, x^j+\de_{i}^j\,1_B\ot 1_A .
\label{per}
\ee
Then, by setting $\varepsilon(x^j)=0$, we get a double $(A,B)$, such that the corresponding operators $\pa^i= Op(x^i)$ are the partial derivatives defined on the polynomial algebra $B$.
 Below, we omit the factors $1_A$ and  $1_B$ in the permutation relations similar to (\ref{per}).

\begin{remark} \rm
In the particular case $N=1$ by slightly modifying  the permutation relations as
$$
y\, x=q\,x\, y+1, \quad q\in \C, \,q\not\in \{0,\pm1\},
$$
we get  the well-known Jackson derivative.
\end{remark}

\noindent
2. Especially, we are interested in a matrix version of the permutation relations (\ref{per}). Consider $N\times N$ matrices $M=\|m_i^j\|$ and $D=\|d_i^j\|$.
Define the algebra $B\ot_{\sigma} A$, where $A=A(D),\, B=B(M)$ (in the brackets we put the generating matrix of the algebra) by the following  system
$$
D_1\, D_2=D_2\, D_1,\quad M_1\, M_2=M_2\, M_1,\quad D_1\,  M_2=M_2\, D_1+P_{12}.
$$

Two first equalities of this system mean that the algebras $A$ and $B$ are commutative. The last equality  (the permutation relations) together with the counit  $\varepsilon(D)=0$ leads to the action
$$
D_1\triangleright M_2=P_{12},\quad \Leftrightarrow \quad \pa_i^j \triangleright m_k^l= \de_i^l\, \de_k^j,\quad{\rm where}\quad \pa_i^j=Op(d_i^j).
$$
The Leibniz rule for the matrix $Op(D)=\|Op(d_i^j)\|$ can be expressed via the coproduct
$$\De(Op(D))=Op(D)\ot I+I\ot Op(D).$$

\noindent
3. Now, consider the double $(A,B)$, with  $A=U(gl(N))$ and $B=T(V)$, where $V$ is  the space  of the covariant representation of the algebra $A$.
Let $\{x_1,...,x_N\}$ be a basis of $V$ and $\{l_i^j\}$ be the corresponding basis of $U(gl(N))$, i.e. such that $l_i^j\triangleright x_k=x_i \de_k^j$.
Then the  relations  between the generators $l_i^j$ can be cast in the following  matrix form
\be
L_1\, L_2-L_2\,L_1=L_1\, P-L_2\, P,\quad L=\|l_i^j\|.
\label{env}
\ee
We impose no relation on the generators $x_i$. The permutation relations are defined as follows
$$
L_2\, \xx_1=\xx_1\, L_2+P_{12}\,\xx_1\quad\Leftrightarrow \quad l_i^j\, x_k=x_k\, l_i^j+ x_i\, \de_k^j.
$$
Hereafter, $\xx$ stands for the column $(x_1,...,x_N)^t$. Note that in this double the algebra $B=T(V)$ can be replaced by one of the  algebras
$\mathrm{Sym}(V)$ or $\mathsf{\Lambda}(V)$. Also, there exist similar doubles with the dual space $V^*$ instead of $V$.

\medskip
\noindent
4. Let $B$ be another copy of the algebra $U(gl(N))$ with a similar basis $m_i^j$ and  $M=\|m_i^j\|$ be the corresponding generating matrix. It meets the system of relations similar
to (\ref{env}). We define two types of the permutation relations by the following formulae
$$
\mathbf{(i):}\,\,L_1\, M_2=M_2\, L_1+M_1\, P_{12}- M_2\, P_{12} \quad \mathrm{or}\quad \mathbf{(ii):}\,\,L_1\, M_2=M_2\, L_1+M_1\, P_{12}.
$$
Then, taking the counit $\varepsilon(L)=0$,  we get the corresponding actions
 $$
\mathbf{(i):}\,\, L_1\triangleright M_2=M_1\, P_{12}-M_2\,P_{12} \quad \mathrm{or}\quad \mathbf{(ii):}\,\, L_1\triangleright M_2=M_1\, P_{12}.
$$

 The above algebra $B=B(M)$ can be replaced by the commutative algebra $\mathrm{Sym}(gl(N))$. Then the relations in the algebra $B$ become $M_1\, M_2=M_2\, M_1$. All other relations
 remain unchanged. In this case the algebra $A=U(gl(N))$ is respectively represented by the adjoint and left vector fields onto  the  algebra   $B=\mathrm{Sym}(gl(N))$.

\medskip
\noindent
 5. The following double was constructed in \cite{3} as a limit case of a double, considered in the next section.  Namely, introduce a double $(A(D),B(N))$, where the generating matrices
$D=\|d_i^j\|$ and $N=\|n_i^j\|$ satisfies the the following systems
$$
D_1\,D_2=D_2\,D_1\,\quad N_1\,N_2-N_2\,N_1 = h\,(N_1\, P_{12}- N_2\,P_{12}),\quad D_1\,N_2=N_2\, D_1+ P_{12}+h \, D_1\,P_{12}.
$$
Thus,  the algebra $A=A(D)$ is commutative and  $B=U(gl(N)_h)$.

The algebra $B\ot_{\sigma} A$ is an NC analog of the HW algebra from the example 2 above. Namely, this algebra is the main ingredient  of our NC $GL(N)$-covariant calculus.

It is convenient to introduce the matrix $\tD=D+h^{-1}\, I$ and simplify the permutation relations to the form:
$$
\tD_1\,N_2=N_2\, \tD_1+ h \, \tD_1\,P_{12}.
$$

By setting $\varepsilon(D)=0$ (and therefore $\varepsilon(\tD)=h^{-1}\, I$),  we get the action of operators  $\pa_i^j=Op(d_i^j)$ on all elements of the algebra $U(gl(N)_h)$.
Note that this action is classical on the generators of the algebra $B(N)$: $\pa_i^j\triangleright n_k^l=\de_i^l\delta_k^j$. Its extension on the higher monomials can be done by means of
the coproduct
$$
\De(\pa_i^j)=\pa_i^j\ot 1+1\ot \pa_i^j-h \sum_k \pa_i^k\ot \pa_k^j.
$$
Observe  that our partial derivatives turn into the usual ones on $\mathrm{Sym}(gl(N))$ as $h\to 0$.

\section{Doubles related to braidings}

1. Let $A=A(T)$, $T=\|t_i^j\| $ and $B=B(M)$, $M=\|M_i^j\|$ be two RTT algebras, corresponding to a Hecke symmetry $R$.
Let us define the permutation relations by the rule
$$
R_{12} T_1 M_2= M_1 T_2\, R_{12} \quad \Leftrightarrow \quad T_1\, M_2= R_{12}^{-1}\,  M_1\, T_2\, R_{12}.
$$
Defining the counit $\varepsilon(T)=I$, we get to the following action
$$
T_1\triangleright M_2=R_{12}^{-1}M_1R_{12}.
$$

\noindent
2. Now, we set $B=T(V)$ and define the permutation relations as follows:
$$
T_1\, \xx_2=R_{12}P_{12} \xx_2\, T_1 \quad \Leftrightarrow \quad t_i^j\, x_k=R_{ik}^{mn }\, x_m\, t_n^j,
$$
where the summation over repeated indices is understood. With the same counit, we have the action
$$
T_1\triangleright \xx_2=R_{12}P_{12}\, \xx_2 \quad \Leftrightarrow \quad t_i^j\triangleright x_k=R_{ik}^{mj }\, x_m.
$$
Computing the action of the elements $t_i^j$ on higher elements from $T(V)$, we arrive to the representations,  which can be constructed via the fusion procedure.
If a Hecke symmetry $R=R(q)\rightarrow P$ at $q\rightarrow 1$, in this limit  we get the trivial representations $t_i^j\to \varepsilon(t_i^j)\, I$ in the both examples above.
Note that the algebra $B=T(V)$ in this construction can be replaced by  $R$-symmetric or $R$-skew-symmetric algebras of the space $V$.

\medskip
\noindent
3. The differential calculus from \cite{6} (see section 7), which is a generalization of the Wess-Zumino calculus on the quantum planes \cite{7},  can be also presented in terms of
a double. Consider a double $(A,B)$ where $B=\mathrm{Sym}_R(V)$ and $A=\mathrm{Sym}_R(V^*)$. Here the space $V^*$ is endowed  with the right dual basis
$\{x^1,...,x^N\}$, i.e. such that $<x_i, x^j>=\de_i^j$. Let us put together  all defining relations of this double:
\be
q\,x_i x_j=R_{ij}^{kl}\,x_k x_l,\quad q\,x^i x^j=R^{ji}_{kl} \,x^l x^k,\quad x^j R_{jk}^{il} \, x_i=h \,\de_k^l+\qq x_k x^l.
\label{Fock}
\ee
Setting $\varepsilon(x^j)=0$, we get $R$-analogs $\pa^i=Op(x^i)$ of the partial derivatives multiplied by  $h$. The above permutation  relations together with the counit
$\varepsilon(x^i)=0$ play the role of  the Leibniz rule for the operators $\pa^i$. The algebra $\mathrm{Sym}_R(V)$ endowed with these operators is an $R$-counterpart
of the bosonic Fock space. In a similar manner an $R$-analog of the fermionic Fock space can be constructed.

\medskip
\noindent
4. Let us consider elements $k_i^j=x_i\, x^j$ and compose the matrix $K=\|k_i^j\|$.

\begin{proposition} In virtue of (\ref{Fock}) the matrix $K$ is subject to the following relation:
\be
R_{12} K_1R_{12} K_1-K_1 R_{12} K_1R_{12}=h\,(R_{12} K_1-K_1 R_{12}).
\label{mRE}
\ee
\end{proposition}
We call the algebra defined by (\ref{mRE}) the {\em modified} RE algebra. If $R$ is an involutive symmetry, the claim above is still valid.
However, only if $R$  is a Hecke symmetry, this algebra is isomorphic to the RE algebra defined by (\ref{RE}). This isomorphism can be defined as follows
\be
L=h\, I-(q-\qq)\, K.
\label{iso}
\ee

Now, consider a double $(A, B)$, where $A=A(K)$ is a modified RE algebra (\ref{mRE}) and $B$ is one of the algebras  $T(V)$, $\mathrm{Sym}_R(V)$,
$\mathsf{\Lambda}_R(V)$. Taking into account (\ref{Fock}) and the identification $k_i^j=x_i\, x^j$ we get the following permutation relations between these algebras
$$
R_{12} K_1 R_{12} \xx_1=\xx_1 K_2+h\, R_{12}\xx_1.
$$
The counit $\varepsilon(K)=0$ leads to  the action
$$
R_{12} K_1 R_{12} \triangleright\xx_1=h\, R_{12} \xx_1.
$$
Assuming that $R=R(q)\to P$ as $q\to 1$, we get the limit action $l_i^j\triangleright x_k=h\, x_i\,\de_k^j$ which coincides with the covariant representation
 of the algebra $U(gl(N)_h)$.

In a similar manner it is possible to define a double with the space $V^*$ instead of $V$ and thus to get the contravariant representation of the modified RE algebra $A$.

In \cite{1} there was described a way of constructing a category of finite dimensional $A$-modules similar to $U(gl(N)_h)$-module. In that construction we used the "braided bi-algebra
structure" of the modified RE algebra and the categorical morphisms transposing the objects $\mathrm{span}(k_i^j)\cong V\ot V^*$ and $M$, where $M$ is an arbitrary object
of the mentioned category.  More precisely,  the corresponding permutation relations are
$$
\sigma(a\ot b)=  (\triangleright_{12}\otimes \mathrm{id})\circ \mathcal{R}_{23}(a_1\otimes a_2\otimes b),\quad {\mathrm{where}}\quad  a_1\otimes a_2 = \Delta(a),
$$
and $\De$ is the coproduct in the modified RE algebra defined in \cite{1}, while $\mathcal{R}$ stands for the braiding (a categorical morphism), transposing
the objects $\mathrm{span}(k_i^j)$ and $M$. Note that in the case related to the quantum group $U_q(sl(N))$, $ \mathcal{R}$ is the product of the usual flip and the image of the
corresponding universal $R$-matrix. However, in general, the mentioned categorical morphism can be constructed via the initial symmetry $R$ without any quantum group.

\medskip
\noindent
5. Let  $A$ be again a modified RE algebra, defined by (\ref{mRE}). The role of $B$ is often attributed to the corresponding RTT algebra. We consider two doubles where
the role of $B=B(M)$, $M=\|m_i^j\|$ is also played by another copy of the RE algebra (in its non-modified form). Define the two types of permutation relations:
\be
\mathbf{(i):}\,\,R_{12} K_1R_{12} M_1=M_1 R_{12} K_1 R_{12}+h\,(R_{12} M_1-R_{12} M_2),
\label{first}
\ee
\be
\mathbf{(ii):}\,\,\ R_{12} K_1R_{12} M_1=M_1 R_{12} K_1 R_{12}+h\,R_{12} M_1.
\label{second}
\ee
The first system of permutation relations defines {\em braided analogs} of the adjoint vector fields. The second one defines {\em braided analogs} of the left vector fields
 (see \cite{2}).

Turn to the double, defined by (\ref{second}). As was shown in \cite{3}, the matrix $D=M^{-1}\, K$ (the matrix $M^{-1}$ can be found via the Cayley-Hamilton identity) and $K$
generate  a double $(A(D),B(M))$, where $D=\|d_i^j\|$ and $M=\|m_i^j\|$,  with the following defining system
$$
R_{12}^{-1} D_1 R_{12}^{-1} D_1=D_1 R_{12}^{-1} D_1R_{12}^{-1},\quad R_{12} M_1 R_{12} M_1=M_1 R_{12} M_1R_{12},
$$
$$
D_1 R_{12} M_1 R_{12}=R_{12}M_1 R_{12}^{-1} D_1+R_{12}.
$$

Now, in this double we replace the  matrix $M$ by $N$ where $M=h\, I-(q-\qq)  N$ and get a double  $(A(D),B(N))$.
The matrix $N=\|n_i^j\|$ generates the modified RE algebra and the permutation relations are as follows
$$
D_1R_{12} N_1R_{12}-R_{12}N_1 R_{12}^{-1} D_1=R_{12}+h\, D_1 R_{12}.
$$
Note that if $R=R(q)$ tends to $P$, in the limit we get the double, exhibited at the end of the previous section.

It is possible to construct similar doubles associated with generalized (braided)  Yangians introduced  in \cite{6}. We plan to consider them elsewhere.

\section{Example of quantization with NC configuration space}

Consider  the last example from section 2 for the case $N=2$ and pass by a change of basis to the algebra $B=U(u(2)_h)$. Making the corresponding change of basis in the algebra
$A$, we get a double $(A,B)$, where $B$ is generated by the elements $t,x,y,z$, subject to the relations
$$
[t,\,x]=[t,\,y]=[t,\,z]=0,\quad [x,\, y]=h\, z,\quad [y,\, z]=h\, x,\quad  [z,\, x]=h\, y.
$$
The commutative algebra $A$ is generated by  $\pa_t$, $\pa_x$, $\pa_y$, $\pa_z$ and the permutation relations read
$$
\begin{array}{l@{\quad}l@{\quad}l@{\quad}l}
[\pa_t,\,t]= {h\over 2}\,\pa_t+1 & [\pa_t,\, x] =-{h\over 2}\,\pa_x &[\pa_t,\, y]=-{h\over 2}\,\pa_y &[\pa_t,\, z]=- {h\over 2}\,\pa_z\\
\rule{0pt}{5mm}
[\pa_x,\, t] = {h\over 2}\,\pa_x &[\pa_x, \,x] = {h\over 2}\,\pa_t+1 & [\pa_x, \, y]={h\over 2}\,\pa_z & [\pa_x, \,z] =- {h\over 2}\,\pa_y \\
\rule{0pt}{5mm}
[\pa_y, \,t] = {h\over 2}\,\pa_y & [\pa_y, \,x] =-{h\over 2}\,\pa_z &[\pa_y, \,y]  = {h\over 2}\,\pa_t+1 & [\pa_y, \,z] = {h\over 2}\,\pa_x\\
\rule{0pt}{5mm}
[\pa_z, \,t] = {h\over 2}\,\pa_z & [\pa_z, \,x] = {h\over 2}\,\pa_y&[\pa_z, \,y] =-{h\over 2}\,\pa_x & [\pa_z, \,z] = {h\over 2}\,\pa_t+1
\end{array}
$$

Introducing the generator $\tpa=\pa_t+2/h \, 1_A$, we can  treat the double $B\ot_{\sigma} A$ as the enveloping algebra of the Lie algebra with 8 generators
$t$, $x$, $y$, $z$, $\tpa$, $\pa_x$, $\pa_y$ and $\pa_z$.

Now, we introduce the so-called {\em quantum radius} by the formula
$$
r_\hh=\sqrt{x^2+y^2+z^2+\hh^2},\quad \hh = h/2i.
$$
In a series of papers we extend\footnote{Note that this extension is not straightforward, since the usual Leibnitz rule for the derivatives on the algebra  $U(u(2)_h)$  is not valid.} the
above partial derivatives onto  any rational functions in $r_\hh$ and some rational functions in $x$, $y$, $z$ and $r_\hh$. This enabled us to extend the procedure of quantization with NC
configurational space onto  a larger algebra. We exhibit two examples of dynamical models obtained by such a quantization. First, consider an equation of Schr\"{o}dinger type
\be
(a\pa_t+b(\pa_x^2+\pa_y^2+\pa_z^2)+\frac{q}{r})\, \psi=0, \label{sch}
\ee
where $a, b, q$ are some constants. In this model  the mentioned quantization is reduced to replacing the usual radius by its quantum counterpart $r_\hh$  and similarly for the derivatives.

Note that in the limit $\hh\to 0$ we get the usual radial Schr\"{o}dinger equation. The problem of finding the eigenvectors and the corresponding eigenvalues for the quantum version of the
Schr\"{o}dinger operator is of interest. In \cite{4} we computed the first correction to the energy of the ground state of (\ref{sch}) caused by the noncommutativity of the space  $U(u(2)_h)$.

Another example is the model of the Dirac monopole in the classical electrodynamics. In a particular case corresponding to the static Dirac monopole the system of Maxwell equations
take the form
$$
\mathrm{rot}\,\mathbf{H}=\mathbf{0}\qquad \mathrm{div}\, \mathbf{ H}=4 g \pi \de(\mathbf{r}),
$$
where $\mathbf{H}=(H_1, H_2, H_3)$ is the magnetic field, $\mathbf{r} = (x,y,z)$, $\mathrm{rot}$ and $\mathrm{div}$ are the curl and divergence respectively.  We succeeded in finding
an NC solution of this system:
$$
\mathbf{H}=\frac{g}{r_\hh(r_\hh^2-\hh^2)}\,\mathbf{r}.
$$
Note that  it tends to the usual Dirac monopole as $\hh\to 0$.

\end{document}